\documentclass[12pt]{article}%
\usepackage{amsmath}
\usepackage{amsfonts}
\usepackage{amssymb}
\usepackage{graphicx}%
\providecommand{\U}[1]{\protect\rule{.1in}{.1in}}
\newtheorem{theorem}{Theorem}

\newtheorem{algorithm}[theorem]{Algorithm}

\newtheorem{case}[theorem]{Case}

\newtheorem{definition}[theorem]{Definition}

\newtheorem{problem}[theorem]{Problem}

\begin{document}

\title{\textbf{Weak and Strong Superiorization: Between Feasibility-Seeking and
Minimization\thanks{Presented at the Tenth Workshop on Mathematical Modelling
of Environmental and Life Sciences Problems, October 16-19, 2014, Constantza,
Romania.\newline http://www.ima.ro/workshop/tenth\_workshop/.}}}
\author{Yair Censor\medskip\\{Department of Mathematics, University of Haifa,}\\Mt. Carmel, Haifa 3498838, Israel}
\date{September 30, 2014. Revised: November 27, 2014.}
\maketitle

\begin{abstract}
We review the superiorization methodology, which can be thought of, in some
cases, as lying between feasibility-seeking and constrained minimization. It
is not quite trying to solve the full fledged constrained minimization
problem; rather, the task is to find a feasible point which is superior (with
respect to an objective function value) to one returned by a
feasibility-seeking only algorithm. We distinguish between two research
directions in the superiorization methodology that nourish from the same
general principle: Weak superiorization and strong superiorization and clarify
their nature.

\end{abstract}

\section{Introduction\label{sect:intro}}

\textbf{What is superiorization}. The superiorization methodology works by
taking an iterative algorithm, investigating its perturbation resilience, and
then, using proactively such permitted perturbations, forcing\ the perturbed
algorithm to do something useful in addition to what it is originally designed
to do. The original unperturbed algorithm is called the \textquotedblleft
Basic Algorithm\textquotedblright\ and the perturbed algorithm is called the
\textquotedblleft Superiorized Version\ of the Basic
Algorithm\textquotedblright.

If the original algorithm\footnote{We use the term \textquotedblleft
algorithm\textquotedblright\ for the iterative processes discussed here, even
for those that do not include any termination criterion. This does not create
any ambiguity because whether we consider an infinite iterative process or an
algorithm with a termination rule is always clear from the context.} is
computationally efficient and useful in terms of the application at hand, and
if the perturbations are simple and not expensive to calculate, then the
advantage of this methodology is that, for essentially the computational cost
of the original Basic Algorithm, we are able to get something more by steering
its iterates according to the perturbations.

This is a very general principle, which has been successfully used in some
important practical applications and awaits to be implemented and tested in
additional fields; see, e.g., the recent papers \cite{rand-conmath,sh14}, for
applications in intensity-modulated radiation therapy and in nondestructive
testing. Although not limited to this case, an important special case of the
superiorization methodology is when the original algorithm is a
feasibility-seeking algorithm, or one that strives to find
constraint-compatible points for a family of constraints, and the
perturbations that are interlaced into the original algorithm aim at reducing
(not necessarily minimizing) a given merit (objective) function. We
distinguish between two research directions in the superiorization methodology
that nourish from the same general principle.

One is the direction when the constraints are assumed to be consistent
(nonempty intersection) and the notion of \textquotedblleft bounded
perturbation resilience\textquotedblright\ is used. In this case one treats
the \textquotedblleft Superiorized Version\ of the Basic
Algorithm\textquotedblright\ as a recursion formula without a stopping rule
that produces an infinite sequence of iterates and asymptotic convergence
questions are in the focus of study.

The second direction does not assume consistency of the constraints but uses
instead a proximity function that measures\ the violation of the constraints.
Instead of seeking asymptotic feasibility, it looks at $\varepsilon
$-compatibility and uses the notion of \textquotedblleft strong perturbation
resilience\textquotedblright. The same core \textquotedblleft Superiorized
Version\ of the Basic Algorithm\textquotedblright\ might be investigated in
each of these directions, but the second is apparently more practical since it
relates better to problems formulated and treated in practice. We use the
terms \textquotedblleft weak superiorization\textquotedblright\ and
\textquotedblleft strong superiorization\textquotedblright\ as a nomenclature
for the first and second directions, respectively\footnote{These terms were
proposed in \cite{cz14-feje}, following a private discussion with our
colleague and coworker in this field Gabor Herman.}.

\textbf{The purpose of this paper}. Since its inception in 2007, the
superiorization method has evolved and gained ground. Quoting and distilling
from earlier publications, we review here the two directions of the
superiorization methodology. A recent review paper on the subject which should
be read together with this paper is Herman's \cite{gth-sup4IA}. Unless
otherwise stated, we restrict ourselves, for simplicity, to the $J$dimensional
Euclidean space $R^{J}$ although some materials below remain valid in Hilbert space.

\textbf{Superiorization} \textbf{related work}. Recent publications on the
superiorization methodology (SM) are devoted to either weak or strong
superiorization, without yet using these terms. They are
\cite{bk13,bdhk07,cdh10,cz12,dhc09,rand-conmath,gh13,hd08,hgdc12,wenma13,ndh12,pscr10}%
, culminating in \cite{sh14} and \cite{cdhst14}. The latter contains a
detailed description of the SM, its motivation, and an up-to-date review of
SM-related previous works scattered in earlier publications, including a
reference to \cite{bdhk07} in which it all started, although without using yet
the terms superiorization and perturbation resilience. \cite{bdhk07} was the
first to propose this approach and implement it in practice, but its roots go
back to \cite{brz06,brz08} where it was shown that if iterates of a
nonexpansive operator converge for any initial point, then its inexact
iterates with summable errors also converge, see also \cite{davidi-thesis}.
Bounded perturbation resilience of a parallel projection method (PPM) was
observed as early as 2001 in \cite[Theorem 2]{combettes2001} (without using
this term). More details on related work appear in \cite[Section 3]{cdhst14}
and in \cite[Section 1]{cz-2-2014}.

\section{The framework\label{sect:frame}}

Let $T$ be a mathematically-formulated problem, of any kind or sort, with
solution set $\Psi_{T}.$ The following cases immediately come to mind although
any $T$ and its $\Psi_{T}$ can potentially be used.

\begin{case}
\label{case:cfp}$T$ is a \textit{convex feasibility problem} (CFP) of the
form: find a vector $x^{\ast}\in\cap_{i=1}^{I}C_{i},$ where $C_{i}\subseteq
R^{J}$ are closed convex subsets of the Euclidean space $R^{J}$. In this case
$\Psi_{T}=\cap_{i=1}^{I}C_{i}.$
\end{case}

\begin{case}
\label{case:2}$T$ is a constrained minimization problem: $\mathrm{minimize}%
\left\{  f(x)\mid x\in\Phi\right\}  $ of an objective function $f$ over a
feasible region $\Phi.$ In this case $\Psi_{T}=\{x^{\ast}\in\Phi\mid
f(x^{\ast})\leq f(x)$ for all $x\in\Phi\}.$
\end{case}

The superiorization methodology is intended for function reduction problems of
the following form.

\begin{problem}
\label{prob:sm}\textbf{The Function Reduction Problem}. Let $\Psi_{T}%
\subseteq{R^{J}}$ be the solution set of some given mathematically-formulated
problem $T$ and let $\phi:R^{J}\rightarrow R$ be an objective function. Let
$\mathcal{A}:R^{J}\rightarrow R^{J}$ be an algorithmic operator that defines
an iterative Basic Algorithm for the solution of $T$. Find a vector $x^{\ast
}\in\Psi_{T}$ whose function $\phi$ value is lesser than that of a point in
$\Psi_{T}$ that would have been reached by applying the Basic Algorithm for
the solution of problem $T.$
\end{problem}

As explained below, the superiorization methodology approaches this problem by
automatically generating from the Basic Algorithm its Superiorized Version.
The so obtained vector $x^{\ast}$ need not be a minimizer of $\phi$ over
$\Psi_{T}.$ Another point to observe is that the very problem formulation
itself depends not only on the data $T,$ $\Psi_{T}$ and $\phi$ but also on the
pair of algorithms -- the original unperturbed Basic Algorithm, represented by
$\mathcal{A},$ for the solution of problem $T,$ and its superiorized version.

A fundamental difference between weak and strong superiorization lies in the
meaning attached to the term \textquotedblleft solution of problem
$T$\textquotedblright\ in Problem \ref{prob:sm}. In weak superiorization
solving the problem $T$ is understood as generating an infinite sequence
$\{x^{k}\}_{k=0}^{\infty}$ that converges to a point $x^{\ast}\in\Psi_{T},$
thus $\Psi_{T}$ must be nonempty. In strong superiorization solving the
problem $T$ is understood as finding a point $x^{\ast}$ that is $\varepsilon
$-compatible with $\Psi_{T},$ for some positive $\varepsilon,$ thus
nonemptiness of $\Psi_{T}$ need not be assumed.

We concentrate in the next sections mainly on Case \ref{case:cfp}.
Superiorization work on Case \ref{case:2}, where $T$ is a maximum likelihood
optimization problem and $\Psi_{T}$ -- its solution set, appears in
\cite{gh13, wenma13, tie}.

\section{Weak superiorization\label{sect:weakSM}}

In weak superiorization the set $\Psi_{T}$ is assumed to be nonempty and one
treats the \textquotedblleft Superiorized Version\ of the Basic
Algorithm\textquotedblright\ as a recursion formula that produces an infinite
sequence of iterates. Convergence questions are studied in their
asymptotically. The SM strives to asymptotically find a point in $\Psi_{T}$
which is superior\textbf{, }i.e., has a lower, but not necessarily minimal,
value of the $\phi$ function, to one returned by the Basic Algorithm that
solves the original problem $T$ only.

This is done by first investigating the bounded perturbation resilience of an
available Basic Algorithm designed to solve efficiently the original problem
$T$ and then proactively using such permitted perturbations to steer its
iterates toward lower values of the $\phi$ objective function while not
loosing the overall convergence to a point in $\Psi_{T}$.

\begin{definition}
\label{def:resilient}\textbf{Bounded perturbation resilience (BPR)}. Let
$\Gamma\subseteq R^{J}$ be a given nonempty set. An algorithmic
operator\ $\mathcal{A}:R^{J}\rightarrow R^{J}$ is said to be \texttt{bounded
perturbations resilient with respect to }$\Gamma$\emph{ }if the following is
true: If a sequence $\{x^{k}\}_{k=0}^{\infty},$ generated by the iterative
process $x^{k+1}=\mathcal{A}(x^{k}),$ for all $k\geq0,$ converges to a point
in $\Gamma$ for all $x^{0}\in R^{J}$, then any sequence $\{y^{k}%
\}_{k=0}^{\infty}$ of points in $R^{J}$ that is generated by $y^{k+1}%
=\mathcal{A}(y^{k}+\beta_{k}v^{k}),$ for all $k\geq0,$ also converges to a
point in $\Gamma$ for all $y^{0}\in R^{J}$ provided that, for all $k\geq0$,
$\beta_{k}v^{k}$ are \texttt{bounded perturbations}, meaning that $\beta
_{k}\geq0$ for all $k\geq0$ such that ${\displaystyle\sum\limits_{k=0}%
^{\infty}}\beta_{k}\,<\infty,$ and that the sequence $\{v^{k}\}_{k=0}^{\infty
}$ is bounded.
\end{definition}

Let $\phi:R^{J}\rightarrow R$ be a real-valued convex continuous function and
let $\partial\phi(z)$ be the subgradient set of $\phi$ at $z$ and, for
simplicity of presentation, assume here that $\Gamma=R^{J}.$ In other specific
cases care must be taken regarding how $\Gamma$ and $\Psi_{T}$ are related.
The following Superiorized Version\ of the Basic Algorithm\textbf{
}$\mathcal{A}$ is based on \cite[Algorithm 4.1]{cz14-feje}.

\begin{algorithm}
\label{alg:super-process}$\left.  {}\right.  $\textbf{Superiorized Version\ of
the Basic Algorithm }$\mathcal{A}$\textbf{.}

\textbf{(0) Initialization}: Let $N$ be a natural number and let $y^{0}\in
R^{J}$ be an arbitrary user-chosen vector.

\textbf{(1)} \textbf{Iterative step}: Given a current iteration vector $y^{k}$
pick an $N_{k}\in\{1,2,\dots,N\}$ and start an inner loop of calculations as follows:

\textbf{(1.1) Inner loop initialization}: Define $y^{k,0}=y^{k}.$

\textbf{(1.2) Inner loop step: }Given $y^{k,n},$ as long as $n<N_{k},$ do as follows:

\textbf{(1.2.1) }Pick a $0<\beta_{k,n}\leq1$ in a way that guarantees that%
\begin{equation}
\sum_{k=0}^{\infty}\sum_{n=0}^{N_{k}-1}\beta_{k,n}<\infty. \label{eq:2.6}%
\end{equation}

\textbf{(1.2.2)} Pick an $\displaystyle s^{k,n}\in\partial\phi(y^{k,n})$ and
define $v^{k,n}$ as follows:%
\begin{equation}
v^{k,n}=\left\{
\begin{array}
[c]{cc}%
-\frac{\displaystyle s^{k,n}}{\displaystyle\left\Vert s^{k,n}\right\Vert }, &
\text{if }0\notin\partial\phi(y^{k,n}),\\
0, & \text{if }0\in\partial\phi(y^{k,n}).
\end{array}
\right.
\end{equation}

\textbf{(1.2.3) }Calculate the perturbed iterate%
\begin{equation}
y^{k,n+1}=y^{k,n}+\beta_{k,n}v^{k,n} \label{eq:2.11}%
\end{equation}
and if $n+1<N_{k}$ set $n\leftarrow n+1$ and go to \textbf{(1.2)}, otherwise
go to \textbf{(1.3)}.

\textbf{(1.3) }Exit the inner loop with the vector $y^{k,N_{k}}$

\textbf{(1.4) }Calculate%
\begin{equation}
y^{k+1}=\mathcal{A}(y^{k,N_{k}}) \label{eq:2.12}%
\end{equation}
set $k\leftarrow k+1$ and go back to \textbf{(1)}.
\end{algorithm}

Let us consider Case \ref{case:cfp} in Section \ref{sect:frame} wherein $T$ is
a convex feasibility problem. The Dynamic String-Averaging Projection (DSAP)
method of \cite{cz12} constitutes a family of algorithmic operators that can
play the role of the above $\mathcal{A}$ in a Basic Algorithm for the solution
of the CFP $T$.

Let $C_{1},C_{2},\dots,C_{m}$ be nonempty closed convex subsets of a Hilbert
space $X$ where $m$ is a natural number. Set $C=\cap_{i=1}^{m}C_{i},$ and
assume $C\neq\emptyset$. For $i=1,2,\dots,m,$ denote by $P_{i}:=P_{C_{i}}$ the
orthogonal (least Euclidean distance) projection onto the set $C_{i}.$ An
index vector is a vector $t=(t_{1},t_{2},\dots,t_{q})$ such that $t_{i}%
\in\{1,2,\dots,m\}$ for all $i=1,2,\dots,q$, whose length is $\ell(t)=q.$ The
product of the individual projections onto the sets whose indices appear in
the index vector $t$ is $P[t]:=P_{t_{q}}\cdots P_{t_{1}}$, called a string operator.

A finite set $\Omega$ of index vectors is called fit if for each
$i\in\{1,2,\dots,m\}$, there exists a vector $t=(t_{1},t_{2},\dots,t_{q}%
)\in\Omega$ such that $t_{s}=i$ for some $s\in\{1,2,\dots,q\}$. Denote by
$\mathcal{M}$ the collection of all pairs $(\Omega,w)$, where $\Omega$ is a
finite fit set of index vectors and $w:\Omega\rightarrow(0,\infty)$ is such
that $\sum_{t\in\Omega}w(t)=1.$

For any $(\Omega,w)\in\mathcal{M}$ define the convex combination of the
end-points of all strings defined by members of $\Omega$%
\begin{equation}
P_{\Omega,w}(x):=\sum_{t\in\Omega}w(t)P[t](x),\;x\in X.
\end{equation}

Let $\Delta\in(0,1/m)$ and an integer $\bar{q}\geq m$ be arbitrary fixed and
denote by $\mathcal{M}_{\ast}\equiv\mathcal{M}_{\ast}(\Delta,\bar{q})$ the set
of all $(\Omega,w)\in\mathcal{M}$ such that the lengths of the strings are
bounded and the weights are all bounded away from zero, i.e.,%
\begin{equation}
\mathcal{M}_{\ast}=\{(\Omega,w)\in\mathcal{M\mid}\text{ }\ell(t)\leq\bar
{q}\text{ and }w(t)\geq\Delta,\text{ }\forall\text{ }t\in\Omega\}.
\end{equation}

\begin{algorithm}
\label{alg:DSAP}$\left.  {}\right.  $\textbf{The DSAP method with variable
strings and variable weights}

\textbf{Initialization}: select an arbitrary $x^{0}\in X$,

\textbf{Iterative step}: given a current iteration vector $x^{k}$ pick a pair
$(\Omega_{k},w_{k})\in\mathcal{M}_{\ast}$ and calculate the next iteration
vector $x^{k+1}$ by%
\begin{equation}
x^{k+1}=P_{\Omega_{k},w_{k}}(x^{k})\text{.}%
\end{equation}

\end{algorithm}

The first prototypical string-averaging algorithmic scheme appeared in
\cite{ceh01} and subsequent work on its realization with various algorithmic
operators includes
\cite{CS08,CS09,ct03,cz-2-2014,crombez,gordon,pen09,pscr10,rhee03}. If in the
DSAP method one uses only a single index vector $t=(1,2,\dots,m)$ that
includes all constraints indices then the fully-sequential Kaczmarz cyclic
projection method is obtained. For linear hyperplanes as constraints sets the
latter is equivalent with the, independently discovered, ART (for Algebraic
Reconstruction Technique) in image reconstruction from projections, see
\cite{GTH}. If, at the other extreme, one uses exactly $m$ one-dimensional
index vectors $t=(i),$ for $i=1,2,\dots,m,$ each consisting of exactly one
constraint index, then the fully-simultaneous projection method of Cimmino is
recovered. In-between these \textquotedblleft extremes\textquotedblright\ the
DSAP method allows for a large arsenal\ of specific feasibility-seeking
projection algorithms. See \cite{bb96, annotated,cccdh10} for more information
on projection methods.

The superiorized version of the DSAP algorithm is obtained by using Algorithm
\ref{alg:DSAP} as the algorithmic operator $\mathcal{A}$ in Algorithm
\ref{alg:super-process}. The following result about its behavior was proved.
Consider the set $C_{min}:=\{x\in C\mid\;\phi(x)\leq\phi(y){\text{ for all }%
}y\in C\},$ and assume that $C_{min}\not =\emptyset.$

\begin{theorem}
\label{thm:fejer}\cite[Theorem 4.1]{cz14-feje} Let $\phi:X\rightarrow R$ be a
convex continuous function, and let $C_{\ast}\subseteq C_{min}$ be a nonempty
subset. Let $r_{0}\in(0,1]$ and $\bar{L}\geq1$ be such that, ${\text{for all
}}x\in C_{\ast}{\text{ and all }}y$ such that$\;||x-y||\leq r_{0},$%
\begin{equation}
|\phi(x)-\phi(y)|\leq\bar{L}||x-y||{,}%
\end{equation}
and suppose that $\{(\Omega_{k},w_{k})\}_{k=0}^{\infty}\subset\mathcal{M}%
_{\ast}.$ Then any sequence $\{y^{k}\}_{k=0}^{\infty},$ generated by the
superiorized version of the DSAP algorithm, converges in the norm of $X$ to a
$y^{\ast}\in C$ and exactly one of the following two alternatives holds:

(a) $y^{\ast}\in C_{min}$;

(b) $y^{\ast}\notin C_{min}$ and there exist a natural number $k_{0}$ and a
$c_{0}\in(0,1)$ such that for each $x\in C_{\ast}$ and each integer $k\geq
k_{0}$,%
\begin{equation}
\Vert y^{k+1}-x\Vert^{2}\leq\Vert y^{k}-x\Vert^{2}-c_{0}\sum_{n=1}^{N_{k}%
-1}\beta_{k,n}.
\end{equation}

\end{theorem}

This shows that $\{y^{k}\}_{k=0}^{\infty}$ is strictly Fej\'{e}r-monotone with
respect to\textbf{ }$C_{\ast},$ i.e., that $\Vert y^{k+1}-x\Vert^{2}<\Vert
y^{k}-x\Vert^{2},$ for all $k\geq k_{0},$ because $c_{0}\sum_{n=1}^{N_{k}%
-1}\beta_{k,n}>0.$ The strict Fej\'{e}r-monotonicity however does not
guarantee convergence to a constrained minimum point but it says that the
so-created feasibility-seeking sequence $\{y^{k}\}_{k=0}^{\infty}$ has the
additional property of getting strictly closer, without necessarily
converging, to the points of a subset of the solution set of of the
constrained minimization problem.

Published experimental results repeatedly confirm that reduction of the value
of the objective function $\phi$ is indeed achieved, without loosing the
convergence toward feasibility, see
\cite{bk13,bdhk07,cdh10,cz12,dhc09,rand-conmath,gh13,hd08,hgdc12,wenma13,ndh12,pscr10}%
. In some of these cases the SM returns a lower value of the objective
function $\phi$ than an exact minimization method with which it is compared,
e.g., \cite[Table 1]{cdhst14}.

\section{Strong superiorization\label{sect:strongSM}}

As in the previous section, let us consider again, Case \ref{case:cfp} in
Section \ref{sect:frame} wherein $T$ is a convex feasibility problem. In this
section we present a restricted version of the SM of \cite{hgdc12} as adapted
to this situation in \cite{cdhst14}. Let $C_{1},C_{2},\dots,C_{m}$ be nonempty
closed convex subsets of $R^{J}$ where $m$ is a natural number and set
$C=\cap_{i=1}^{m}C_{i}$. We do not assume that $C\neq\emptyset,$ but only that
there is some nonempty subset $\Lambda\in R^{J}$ such that $C\subseteq
\Lambda.$ Instead of the nonemptiness assumption we associate with the family
of constraints $\{C_{i}\}_{i=1}^{m}$ a \textit{proximity function} ${Prox}%
_{C}:\Lambda\rightarrow\mathbb{R}_{+}$ that is an indicator of how
incompatible an $x\in\Lambda$ is with the constraints. For any given
$\varepsilon>0$, a point $x\in\Lambda$ for which ${Prox}_{C}(x)\leq
\varepsilon$ is called an $\varepsilon$\textit{-compatible solution} for $C$.
We further assume that we have a feasibility-seeking \textit{algorithmic
operator} $\mathcal{A}:R^{J}\rightarrow\Lambda$, with which we define the
Basic Algorithm as the iterative process%
\begin{equation}
x^{k+1}=\mathcal{A}(x^{k}),\text{ for all }k\geq0,\text{ for an arbitrary
}x^{0}\in\Lambda.
\end{equation}
The following definition helps to evaluate the output of the Basic Algorithm
upon termination by a stopping rule. This definition as well as most of the
remainder of this section appeared in \cite{hgdc12}.

\begin{definition}
\textbf{The }$\varepsilon$\textbf{-output of a sequence. }Given $C\subseteq
\Lambda\subseteq R^{J}$, a proximity function ${Prox}_{C}:\Lambda\rightarrow
R_{+}$, a sequence $\left\{  x^{k}\right\}  _{k=0}^{\infty}\subset\Lambda$ and
an $\varepsilon>0,$ then an element $x^{K}$ of the sequence which has the
properties: (i) ${Prox}_{C}\left(  x^{K}\right)  \leq\varepsilon,$ and (ii)
${Prox}_{C}\left(  x^{k}\right)  >\varepsilon$ for all $0\leq k<K,$ is called
an $\varepsilon$\texttt{-output of the sequence }$\left\{  x^{k}\right\}
_{k=0}^{\infty}$\texttt{ with respect to the pair }$(C,$\texttt{ }${Prox}%
_{C})$.
\end{definition}

We denote the $\varepsilon$-output by $O\left(  C,\varepsilon,\left\{
x^{k}\right\}  _{k=0}^{\infty}\right)  =x^{K}.$ Clearly, an $\varepsilon
$-output $O\left(  C,\varepsilon,\left\{  x^{k}\right\}  _{k=0}^{\infty
}\right)  $ of a sequence $\left\{  x^{k}\right\}  _{k=0}^{\infty}$ might or
might not exist, but if it does, then it is unique. If $\left\{
x^{k}\right\}  _{k=0}^{\infty}$ is produced by an algorithm intended for the
feasible set $C,$ such as the Basic Algorithm, without a termination
criterion, then $O\left(  C,\varepsilon,\left\{  x^{k}\right\}  _{k=0}%
^{\infty}\right)  $ is the \textit{output} produced by that algorithm when it
includes the termination rule to stop when an $\varepsilon$-compatible
solution for $C$ is reached.

\begin{definition}
\textbf{Strong perturbation resilience. }Assume that we are given a
$C\subseteq\Lambda$, a proximity function ${Prox}_{C}$, an algorithmic
operator $\mathcal{A}$ and an $x^{0}\in\Lambda$. We use $\left\{
x^{k}\right\}  _{k=0}^{\infty}$ to denote the sequence generated by the Basic
Algorithm when it is initialized by $x^{0}$. The Basic Algorithm is said to
be\texttt{ strongly perturbation resilient} iff the following hold: (i) there
exist an $\varepsilon>0$ such that the $\varepsilon$-output $O\left(
C,\varepsilon,\left\{  x^{k}\right\}  _{k=0}^{\infty}\right)  $ exists for
every $x^{0}\in\Lambda$; (ii) for every $\varepsilon>0,$ for which the
$\varepsilon$-output $O\left(  C,\varepsilon,\left\{  x^{k}\right\}
_{k=0}^{\infty}\right)  $ exists for every $x^{0}\in\Lambda$, we have also
that the $\varepsilon^{\prime}$-output $O\left(  C,\varepsilon^{\prime
},\left\{  y^{k}\right\}  _{k=0}^{\infty}\right)  $ exists for every
$\varepsilon^{\prime}>\varepsilon$ and for every sequence $\left\{
y^{k}\right\}  _{k=0}^{\infty}$ generated by%
\begin{equation}
y^{k+1}=\mathcal{A}\left(  y^{k}+\beta_{k}v^{k}\right)  ,\text{ for all }%
k\geq0, \label{eq:perturb}%
\end{equation}
where the vector sequence $\left\{  v^{k}\right\}  _{k=0}^{\infty}$ is bounded
and the scalars $\left\{  \beta_{k}\right\}  _{k=0}^{\infty}$ are such that
$\beta_{k}\geq0$, for all $k\geq0,$ and $\sum_{k=0}^{\infty}\beta_{k}<\infty$.
\end{definition}

A theorem which gives sufficient conditions for strong perturbation resilience
of the Basic Algorithm has been proved in \cite[Theorem 1]{hgdc12}. Along with
the $C\subseteq R^{J}$, we look at the objective function $\phi:R^{J}%
\rightarrow R$, with the convention that a point in $R^{J}$ for which the
value of $\phi$ is smaller is considered \textit{superior} to a point in
$R^{J}$ for which the value of $\phi$ is larger. The essential idea of the SM
is to make use of the perturbations of (\ref{eq:perturb}) to transform a
strongly perturbation resilient Basic Algorithm that seeks a
constraints-compatible solution for $C$ into its Superiorized Version whose
outputs are equally good from the point of view of constraints-compatibility,
but are superior (not necessarily optimal) according to the objective function
$\phi$.

\begin{definition}
\label{def:nonascend}Given a function $\phi:R^{J}\rightarrow R$ and a point
$y\in R^{J}$, we say that a vector $d\in R^{J}$ is \texttt{nonascending}
\texttt{for }$\phi$\texttt{ at }$y$ iff $\left\Vert d\right\Vert \leq1$ and
there is a $\delta>0$ such that for all $\lambda\in\left[  0,\delta\right]  $
we have $\phi\left(  y+\lambda d\right)  \leq\phi\left(  y\right)  .$
\end{definition}

Obviously, the zero vector is always such a vector, but for superiorization to
work we need a sharp inequality to occur in (\ref{def:nonascend}) frequently
enough. The Superiorized Version of the Basic Algorithm assumes that we have
available a summable sequence $\left\{  \eta_{\ell}\right\}  _{\ell=0}%
^{\infty}$ of positive real numbers (for example, $\eta_{\ell}=a^{\ell}$,
where $0<a<1$) and it generates, simultaneously with the sequence $\left\{
y^{k}\right\}  _{k=0}^{\infty}$ in $\Lambda$, sequences $\left\{
v^{k}\right\}  _{k=0}^{\infty}$ and $\left\{  \beta_{k}\right\}
_{k=0}^{\infty}$. The latter is generated as a subsequence of $\left\{
\eta_{\ell}\right\}  _{\ell=0}^{\infty}$, resulting in a nonnegative summable
sequence $\left\{  \beta_{k}\right\}  _{k=0}^{\infty}$. The algorithm further
depends on a specified initial point $y^{0}\in\Lambda$ and on a positive
integer $N$. It makes use of a logical variable called \textit{loop}\emph{.
}The Superiorized Version of the Basic Algorithm is presented next by its pseudo-code.

\begin{algorithm}
\label{alg_super}\textbf{Superiorized Version of the Basic Algorithm}
\end{algorithm}

\begin{enumerate}
\item \textbf{set} $k=0$

\item \textbf{set} $y^{k}=y^{0}$

\item \textbf{set} $\ell=-1$

\item \textbf{repeat}

\item $\qquad$\textbf{set} $n=0$

\item $\qquad$\textbf{set} $y^{k,n}=y^{k}$

\item $\qquad$\textbf{while }$n$\textbf{$<$}$N$

\item $\qquad$\textbf{$\qquad$set }$v^{k,n}$\textbf{ }to be a nonascending
vector for $\phi$ at $y^{k,n}$

\item $\qquad$\textbf{$\qquad$set} \emph{loop=true}

\item $\qquad$\textbf{$\qquad$while}\emph{ loop}

\item $\qquad\qquad\qquad$\textbf{set $\ell=\ell+1$}

\item $\qquad\qquad\qquad$\textbf{set} $\beta_{k,n}=\eta_{\ell}$

\item $\qquad\qquad\qquad$\textbf{set} $z=y^{k,n}+\beta_{k,n}v^{k,n}$

\item $\qquad\qquad\qquad$\textbf{if }$\phi\left(  z\right)  $\textbf{$\leq$%
}$\phi\left(  y^{k}\right)  $\textbf{ then }

\item $\qquad\qquad\qquad\qquad$\textbf{set }$n$\textbf{$=$}$n+1$

\item $\qquad\qquad\qquad\qquad$\textbf{set }$y^{k,n}$\textbf{$=$}$z$

\item $\qquad\qquad\qquad\qquad$\textbf{set }\emph{loop = false}

\item $\qquad$\textbf{set }$y^{k+1}$\textbf{$=$}$\mathcal{A}\left(
y^{k,N}\right)  $

\item $\qquad$\textbf{set }$k=k+1$
\end{enumerate}

\begin{theorem}
\label{theorem4.5-1} Any sequence $\left\{  y^{k}\right\}  _{k=0}^{\infty}$,
generated by the Superiorized Version of the Basic Algorithm, Algorithm
\ref{alg_super}, satisfies (\ref{eq:perturb}). Further, if, for a given
$\varepsilon>0,$ the $\varepsilon$-output $O\left(  C,\varepsilon,\left\{
x^{k}\right\}  _{k=0}^{\infty}\right)  $ of the Basic Algorithm exists for
every $x^{0}\in\Lambda$, then every sequence $\left\{  y^{k}\right\}
_{k=0}^{\infty}$, generated by the Algorithm \ref{alg_super}, has an
$\varepsilon^{\prime}$-output $O\left(  C,\varepsilon^{\prime},\left\{
y^{k}\right\}  _{k=0}^{\infty}\right)  $ for every $\varepsilon^{\prime
}>\varepsilon$.
\end{theorem}

The proof of this theorem follows from the analysis of the behavior of the
Superiorized Version of the Basic Algorithm in \cite[pp. 5537--5538]{hgdc12}.
In other words, Algorithm \ref{alg_super} produces outputs that are
essentially as constraints-compatible as those produced by the original Basic
Algorithm. However, due to the repeated steering of the process by lines 7 to
17 toward reducing the value of the objective function $\phi$, we can expect
that its output will be superior (from the point of view of $\phi$) to the
output of the (unperturbed) Basic Algorithm.

Algorithms \ref{alg:super-process} and \ref{alg_super} are not identical. For
example, the first employes negative subgradients while the second allows to
use any nonascending directions of $\phi.$ Nevertheless, they are based on the
same leading principle of the superiorization methodology. Comments on the
differences between them can be found in \cite[Remark 4.1]{cz14-feje}. While
experimental work has repeatedly demonstrated benefits of the SM, the Theorems
\ref{thm:fejer} and \ref{theorem4.5-1} related to these superiorized versions
of the Basic Algorithm, respectively, leave much to be desired in terms of
rigorously analyzing the behavior of the SM under various conditions.

\section{Concluding comments\label{sect: conclusion}}

In many mathematical formulations of significant real-world technological or
physical problems, the objective function is exogenous to the modeling process
which defines the constraints. In such cases, the faith\ of the modeler in the
usefulness of an objective function for the application at hand is limited
and, as a consequence, it is probably not worthwhile to invest too much
resources in trying to reach an exact constrained minimum point. This is an
argument in favor of using the superiorization methodology for practical
applications. In doing so the amount of computational efforts invested
alternatingly between performing perturbations and applying the Basic
Algorithm's algorithmic operator can, and needs to, be carefully controlled in
order to allow both activities to properly influence\ the outcome. Better
theoretical insights into the behavior of weak and of strong superiorization
as well as better ways of implementing the methodology are needed and await to
be developed.

\end{document}